\documentclass{amsart}
\usepackage{amsmath,amssymb,amsfonts,amsbsy}

\newtheorem{thm}{Theorem}[section]

\newtheorem{example}[thm]{Example}
\newtheorem{lem}[thm]{Lemma}
\newtheorem{prop}[thm]{Proposition}
\newtheorem{defn}[thm]{Definition}

\def\gfrac#1#2{\left[\begin{array}{c}#1\\#2\end{array}\right]}

\newcommand{\bi}{{\mathbf i}}
\newcommand{\bX}{{\mathbf X}}
\newcommand{\bY}{{\mathbf Y}}
\newcommand{\bUpsilon}{{\boldsymbol\Upsilon}}
\newcommand{\bPhi}{{\boldsymbol\Phi}}

\newcommand{\cG}{{\mathcal G}}
\newcommand{\cH}{{\mathcal H}}

\DeclareMathOperator{\dlog}{dlog}
\DeclareMathOperator{\supp}{supp}

\DeclareMathOperator{\mH}{H}
\DeclareMathOperator{\res}{res}

\begin{document}

\title{Parameters Changes for Generalized Power Series}
\author{I-Chiau Huang}
\address{Institute of Mathematics, Academia Sinica, Nankang, Taipei 11529, Taiwan, R.O.C.}
\email{ichuang@math.sinica.edu.tw}

\begin{abstract}
Differentials are introduced to the method of generating functions for generalized power series
with exponents in a totally ordered Abelian group. Logarithmic analogue of cohomology residues
is defined to equate coefficients.
\end{abstract}


\maketitle


\section{Introduction}


The goal of this paper is to provide an algebraic framework for a combinatorial
phenomenon arising from a resemblance of variables changes. It is motivated by a well-known
formula of Jacobi~\cite{jac:rasi} stated and generalized in this paper as
Theorem~\ref{thm:Jacobi} in its modern guise. Our central idea is to find a new notion of
differential and a generalization of variable, with which Jacobians
appear naturally. To obtain combinatorial information from this algebraic framework, we define an
analogue of cohomology residue map \cite{hu:arci}. The new residue map also fits our philosophy
that a residue comes from a differential. Another recent interpretation of 
Jacobi's formula can be found in \cite{xin:rtmns}, where differentials are lacking.

To gain a perspective of this paper, it is helpful to look at the method of generating
functions whose foundation is built up by rings of formal power series and the 
operation of equating coefficients. While elegant and easy to implement, the effect of 
variable changes is not clear in the method without the notion of  K\"{a}hler differentials. 
For instance, the role of the Jacobian occurring in a variables change is not transparent. 
The situation can be improved by meromorphic differentials, with which contour integrations give an 
alternative way to take coefficients \cite{eg:irccs}. In the analytic procedure of the method
of generating functions, extra attentions are paid to conditions without combinatorial
significance such as convergence of sequences and paths for integrations.

Removing unnecessary analytic constraints, the author arrives at certain cohomology classes of
separated differentials \cite{hu:arci}. The process of integration is replaced by cohomology
residue maps, which play a significant role in Grothendieck duality theory. Our new algebraic
framework is supported by the fact that formal power series rings are in fact a notion free
from variables. From this point of view, Lagrange inversion formulae are simply a phenomenon
of variables changes \cite{hu:rpsr}; and pairs of inverse relations are a phenomenon of
Schauder bases changes~\cite{hu:irsb}. The cohomology residues come in as an amenable
tool for realizing these phenomena. 

The formalism of cohomology residues is simple. For instance,
$$
\res\gfrac{\Phi d X_1\wedge\cdots\wedge d X_n}{X_1,\cdots,X_n}=
\text{ the constant term of $\Phi$}
$$
for a power series $\Phi\in\kappa[[X_1,\cdots,X_n]]$ with coefficients in a field $\kappa$. 
Although working well on wide range of problems in combinatorial analysis
\cite{hu:rmca}, a variable change from $X_i$ to $X_i^{-1}$ is not available. In this paper, we
work on a field $\kappa[[e^\cG]]$ of generalized power series with exponents in a totally ordered 
Abelian group $\cG$ and coefficients in $\kappa$ (see Section~\ref{sec:gps} for a review). 
The notion of variables is extended to include their inverses. The logarithmic analogue 
$$
\res\gfrac{\Phi \dlog X_1\wedge\cdots\wedge \dlog X_n}{\log X_1,\cdots,\log X_n}
$$
of residues is defined, even for a field of positive characteristic. The framework consists
of differentials (Section~\ref{sec:diff}), parameters, generalized fractions
and residue maps (Section~\ref{sec:pr}). The useful Jacobi's formula (Theorem~\ref{thm:Jacobi})
and Dyson's conjecture (Section~\ref{sec:appl})
are in fact a phenomenon of parameters changes. These interpretations are seen naturally in our
framework.


\section{Generalized power series}\label{sec:gps}


We recall the definition and basic properties of generalized power series. For details of  
proofs, the reader is referred to \cite[Chapter 13, \S 2]{pas:asgr}. 
Generalized power series are called Malcev-Neumann series in \cite{xin:rtmns}. See
\cite{ber:cnrhlgr,ber:hacnrhlgr} for historical remarks on choices of these names.

Let $\cG$ be an Abelian group
and $\kappa$ be a field, whose elements are called scalars. The set
$$
\kappa\{e^\cG\}:=\{\sum_{g\in\cG}a_ge^g\colon a_g\in\kappa\}
$$
of formal sums is a $\kappa$-vector space with termwise addition and multiplication
\begin{eqnarray}\label{eq:addition}
\sum_{g\in\cG}a_ge^g+\sum_{g\in\cG}b_ge^g&:=&\sum_{g\in\cG}(a_g+b_g)e^g,\\
b\sum_{g\in\cG}a_ge^g&:=&\sum_{g\in\cG}(ba_g)e^g\nonumber.
\end{eqnarray}
For $\Phi=\sum_{g\in\cG}a_ge^g\in\kappa\{e^\cG\}$, we call $a_g$ the
$\kappa$-coefficient of $\Phi$ at $e^g$. The $\kappa$-coefficient of $\Phi$ at $1$ is called
the constant term of $\Phi$ in $\kappa$. An element of $\kappa\{e^\cG\}$ is determined by its
$\kappa$-coefficients. For an element of the form $Y=e^g$ with $g\in\cG$, we use the notation
$$
\log Y:=g.
$$
The {\em support} of $\Phi$ is defined as
$$
\supp\Phi:=\{g\in\cG\colon\text
{ the $\kappa$-coefficient of $\Phi$ at $e^g$ is not zero}\}.
$$ 
A totally ordered Abelian group is an Abelian group together with a total order compatible
with the group structure.
\begin{defn}[generalized power series]
Let $\cG$ be a totally ordered Abelian group. We define
$$
\kappa[[e^\cG]]:=
\{\Phi\in\kappa\{e^\cG\}\colon\text{$\supp\Phi$ is well-ordered}\}.
$$
An element in $\kappa[[e^\cG]]$ is called a generalized power series
with exponents in $\cG$ and coefficients in $\kappa$. 
\end{defn}
Recall that a subset $A$ of $\cG$ is well-ordered if every non-empty subset of $A$ has a
smallest element. 

\begin{lem}\label{lem:9413451}
Let $I_1,\cdots,I_n$ be well-ordered subsets of $\cG$ and $g\in\cG$. The equation 
$x_1+\cdots+x_n=g$ has finitely many solutions $(x_1,\cdots,x_n)$ with $x_i\in I_i$. 
The set $I_1+\cdots+I_n=\{a_1+\cdots+a_n\,|\,a_i\in I_i\}$ is well-ordered. 
\end{lem}

So we may define multiplication
$$
(\sum_{g\in\cG}a_ge^g)(\sum_{g\in\cG}b_ge^g):=
\sum_{g\in\cG}\left(\sum_{g_1+g_2=g}a_{g_1}b_{g_2}\right)e^g
$$
for generalized power series. Together with the addition~(\ref{eq:addition}), $\kappa[[e^\cG]]$
is a commutative ring with the unit $1:=e^0$. A generalized power series $\sum a_ge^g$ is positive
if $a_g=0$ for all $g\leq 0$. 

Let $\Phi=\sum a_ge^g$ be a generalized power series. There are no strictly
decreasing infinite sequences in $\supp\Phi$. If $\Phi$ is positive, for a fixed $g\in\cG$, there are only 
finitely many $i$ such that
$$
\sum_{g_1+\cdots+g_i=g}a_{g_1}\cdots a_{g_i}\neq 0.
$$ 
Given scalars $c_i$, we can define a generalized power series $c_0+c_1\Phi+c_2\Phi^2+\cdots$, 
whose $\kappa$-coefficient at $e^g$ is 
$$
\sum_i\left(c_i\sum_{g_1+\cdots+g_i=g}a_{g_1}\cdots a_{g_i}\right).
$$
If the characteristic of $\kappa$ is zero, we define
$$
\log(1+\Phi):=\sum_{\ell=1}^\infty(-1)^{\ell+1}\frac{\Phi^\ell}{\ell}.
$$

A non-zero generalized power series $\Psi$ can be factorized uniquely as
$$
\Psi=aY(1+\tilde{\Psi}),
$$
where $a$ is a scalar, $\tilde{\Psi}$ is a positive generalized power series and $Y=e^g$ for some $g\in\cG$.
Indeed, $g$ is the smallest element of $\supp\Psi$, $a$ is the $\kappa$-coefficient
of $\Psi$ at $e^g$ and $\tilde{\Psi}=a^{-1}e^{-g}\Psi-1$. We call $a$ the leading 
$\kappa$-coefficient of $\Psi$. In this paper, the factorization of a non-zero 
generalized power series refers to the representation of the above form. $\Psi$ is 
invertible with the inverse 
$a^{-1}e^{-g}(1-\tilde{\Psi}+\tilde{\Psi}^2-\cdots)$. We conclude that
$\kappa[[e^\cG]]$ is a field.

\begin{example}[Laurent series]
The field $\kappa[[e^{\mathbb Z}]]$ with the usual order on $\mathbb Z$ is isomorphic to the
field $\kappa((X))$ of Laurent series.
\end{example}
\begin{example}[Hahn \mbox{\cite{hah:nags}, see \cite[p. 445-499]{hah:ga1}}]
The field $\kappa[[e^{\mathbb Q}]]$ with the usual order on $\mathbb Q$ is of
particular interest to algebraic geometers, since it  contains an algebraic closure of
$\kappa((X))$ if $\kappa$ is algebraically closed.
\end{example}
\begin{example}[iterated Laurent series]\label{example:z2}
Let $\cG_i={\mathbb Z}\oplus{\mathbb Z}$ ($i=1,2$) and let $X=e^{(1,0)}$ and $Y=e^{(0,1)}$ be 
elements of $\kappa\{e^{{\mathbb Z}\oplus{\mathbb Z}}\}$. With the order
\begin{eqnarray*}
(m_1,n_1)\geq(m_2,n_2)&\Longleftrightarrow&
m_1>m_2 \text{ or } m_1=m_2\text{ $\&$ } n_1\geq n_2
\end{eqnarray*} 
on $\cG_1$, $\kappa[[e^{\cG_1}]]$ is isomorphic to the field $\kappa((Y))((X))$ of iterated
Laurent series. With the order
\begin{eqnarray*}
(m_1,n_1)\geq(m_2,n_2)&\Longleftrightarrow&
n_1>n_2 \text{ or } n_1=n_2\text{ $\&$ } m_1\geq m_2
\end{eqnarray*} 
on $\cG_2$, $\kappa[[e^{\cG_2}]]$ is isomorphic to the field $\kappa((X))((Y))$ of iterated
Laurent series. As subsets of $\kappa\{e^{{\mathbb Z}\oplus{\mathbb Z}}\}$, $\kappa[[e^{\cG_1}]]$ 
and $\kappa[[e^{\cG_2}]]$ are different. The inverse of $X+Y$ in $\kappa[[e^{\cG_1}]]$ is
$$
Y^{-1}(1-XY^{-1}+X^2Y^{-2}-X^3Y^{-3}+\cdots)
$$
and that in $\kappa[[e^{\cG_2}]]$ is
$$
X^{-1}(1-X^{-1}Y+X^{-2}Y^2-X^{-3}Y^3+\cdots).
$$
\end{example}

Let $\cH$ be a subgroup of $\cG$ with the induced order. In the rest of this
paper, we assume that the quotient group $\cG/\cH$ is a free Abelian group of rank $n$. In other
words, there exist $u_1,\cdots,u_n\in\cG$ such that every element in $\cG$ can be written
uniquely as $h+s_1u_1+\cdots+s_nu_n$ with $h\in\cH$ and $s_i\in\mathbb Z$. We say also that
$\cG$ is
generated freely by $\cH$ and $u_1,\cdots,u_n$.
\begin{defn}[variable]
$e^{u_1},\cdots,e^{u_n}$ are variables of $\kappa[[e^\cG]]$ over $\kappa[[e^\cH]]$ if $\cG$ is
generated freely by $\cH$ and $u_1,\cdots,u_n$.
\end{defn}
The cardinalities of any sets of variables of $\kappa[[e^\cG]]$ over $\kappa[[e^\cH]]$ are the 
same. We often use the notation $X_i=e^{u_i}$. We say also that
$\kappa[[e^\cG]]$ is generated by $\kappa[[e^\cH]]$ and the variables $X_1,\cdots,X_n$. For an element
$$
\Psi=\sum_{\overset{h\in\cH}{j_1,\cdots,j_n\in\mathbb Z}}
a_{h,j_1,\cdots,j_n}e^hX_1^{j_1}\cdots X_n^{j_n}
$$
in $\kappa[[e^\cG]]$ and fixed $j_1,\cdots,j_n\in\mathbb Z$, we call 
$\varphi_{j_1,\cdots,j_n}:=\sum_{h\in\cH}a_{h,j_1,\cdots,j_n}e^h$
the $\kappa[[e^\cH]]$-coefficient of $\Psi$ at the monomial $X_1^{j_1}\cdots X_n^{j_n}$. The
$\kappa[[e^\cH]]$-coefficient of $\Psi$ at $1$ (that is, $\varphi_{0,\cdots,0}$) is independent of the choice of variables and is
called the constant term of $\Psi$ in $\kappa[[e^\cH]]$. Indeed, fo any $h\in\cH$, the
$\kappa$-coefficient of $\varphi_{0,\cdots,0}$ at $e^h$ equals to that of $\Psi$ at $e^h$. For a
given set of variables,
$\Psi$ is determined by its $\kappa[[e^\cH]]$-coefficients. We use the notation
$\Psi=\Psi(X_1,\cdots,X_n)$ to indicate that it is represented as the form
$$
\Psi=\sum_{j_1,\cdots,j_n\in\mathbb Z}\varphi_{j_1,\cdots,j_n}X_1^{j_1}\cdots X_n^{j_n}
$$
with $\varphi_{j_1,\cdots,j_n}\in\kappa[[e^\cH]]$.


\section{Differentials}\label{sec:diff}

 
A derivation on $\kappa[[e^\cG]]$ is a map 
$D$ from $\kappa[[e^\cG]]$ to a $\kappa[[e^\cG]]$-vector space which satisfies 
$D(\Phi_1+\Phi_2)=D(\Phi_1)+D(\Phi_2)$ and 
$D(\Phi_1\Phi_2)=\Phi_1D(\Phi_2)+\Phi_2D(\Phi_1)$ for all $\Phi_1,\Phi_2\in\kappa[[e^\cG]]$. 
Recall that $\cH$ is a subgroup of $\cG$. A derivation $D$ on $\kappa[[e^\cG]]$ 
is a $\kappa[[e^\cH]]$-derivation if $D(\varphi)=0$ for all $\varphi\in\kappa[[e^\cH]]$.

\begin{defn}[partial derivation]
Let $X_1,\cdots,X_n$ be a set of variables of $\kappa[[e^\cG]]$ over $\kappa[[e^\cH]]$. The 
partial derivation on $\kappa[[e^\cG]]$ with respect to $X_i$ is the well-defined $\kappa[[e^\cH]]$-derivation
$$
\frac{\partial}{\partial X_i}\colon\kappa[[e^\cG]]\to\kappa[[e^\cG]]
$$
given by
$$
\sum\varphi_{j_1,\cdots,j_n}X_1^{j_1}\cdots X_n^{j_n}
\mapsto
\sum j_i\varphi_{j_1,\cdots,j_n}
X_1^{j_1}\cdots X_{i-1}^{j_{i-1}}X_i^{j_i-1}X_{i+1}^{j_{i+1}}\cdots X_n^{j_n}.
$$
\end{defn}
\begin{defn}[compatibility with partial derivations]
Let $Y_1,\cdots,Y_n$ be a set of variables of $\kappa[[e^\cG]]$ over $\kappa[[e^\cH]]$.
A $\kappa[[e^\cH]]$-derivation $D$ on $\kappa[[e^\cG]]$ is compatible with partial derivations 
$\partial/\partial Y_1,\cdots,\partial/\partial Y_n$,
if 
\begin{equation}\label{eq:part}
D(\Phi)=\sum_{i=1}^n\frac{\partial\Phi}{\partial Y_i}D(Y_i)
\end{equation}
for any $\Phi\in\kappa[[e^\cG]]$. $D$ is compatible with partial derivations if it is compatible with 
$\partial/\partial Y_1,\cdots,\partial/\partial Y_n$ for any set of variables $Y_1,\cdots,Y_n$. 
\end{defn}
A $\kappa[[e^\cH]]$-derivation compatible with partial derivations 
$\partial/\partial Y_1,\cdots,\partial/\partial Y_n$ is determined by its values at
$Y_1,\cdots,Y_n$.
\begin{lem}
A partial derivation is compatible with partial derivations.
\end{lem}
\begin{proof}
Let $X_1,\cdots,X_n$ and $Y_1,\cdots,Y_n$ 
be two sets of variables of $\kappa[[e^\cG]]$ over $\kappa[[e^\cH]]$ with the relation 
$$
\begin{cases}
Y_i=e^{h_i}X_1^{s_{i1}}\cdots X_n^{s_{in}},\\
X_i=e^{h'_i}Y_1^{t_{i1}}\cdots Y_n^{t_{in}},
\end{cases}
$$
where $h_i,h_i'\in\cH$. Note that the matrix $(t_{ij})$ is the inverse of $(s_{ij})$. To 
show that $\partial/\partial X_j$ is compatible with partial derivations, 
we check first the relation (\ref{eq:part}) for $\Phi=X_k$:
\begin{eqnarray*}
\sum_{i=1}^n\frac{\partial X_k}{\partial Y_i}\frac{\partial Y_i}{\partial X_j}&=&
\sum_{i=1}^n\frac{e^{h'_k}\partial(Y_1^{t_{k1}}\cdots Y_n^{t_{kn}})}{\partial Y_i}
\frac{e^{h_i}\partial (X_1^{s_{i1}}\cdots X_n^{s_{in}})}{\partial X_j}\\
&=&
\sum_{i=1}^nt_{ki}s_{ij}\frac{X_k}{X_j}=\delta_{kj}=\frac{\partial X_k}{\partial X_j}.
\end{eqnarray*}
For the general case, it suffices to prove that the coefficients of both sides of 
(\ref{eq:part}) at any fixed $g\in\cG$ are the same. Since the coefficients involve 
only finitely many $\kappa[[e^\cH]]$-coefficients at monomials in $X_1,\cdots,X_n$, 
the general case is reduced to the special case that $\Phi$ equals to a finite sum of 
elements of the form $\varphi X_1^{k_1}\cdots X_n^{k_n}$ with $\varphi\in\kappa[[e^\cH]]$.
From the defining properties of derivations, the special case is further reduced to 
the case $\Phi=X_k$ that we just proved.
\end{proof}
\begin{prop}[criterion of compatibility]
Let $D$ be a $\kappa[[e^\cH]]$-derivation on $\kappa[[e^\cG]]$. If $D$ is compatible with
$\partial/\partial X_1,\cdots,\partial/\partial X_n$
for one set of variables $X_1,\cdots,X_n$, then $D$ is compatible with partial derivations.
\end{prop}
\begin{proof}
Let $Y_1,\cdots,Y_n$ 
be another set of variables of $\kappa[[e^\cG]]$ over $\kappa[[e^\cH]]$. The proposition follows
from the straightforward computations:
$$
\sum_{i=1}^n\frac{\partial\Phi}{\partial Y_i}D(Y_i)=
\sum_{i=1}^n\sum_{j=1}^n\frac{\partial\Phi}{\partial Y_i}
\frac{\partial Y_i}{\partial X_j}D(X_j)=
\sum_{j=1}^n\frac{\partial\Phi}{\partial X_j}D(X_j)=D(\Phi).
$$
\end{proof}

\begin{defn}[differentials]
A $\kappa[[e^\cG]]$-vector space $\Omega_{\cG/\cH}$ together with a $\kappa[[e^\cH]]$-derivation
$d\colon\kappa[[e^\cG]]\to\Omega_{\cG/\cH}$ is the vector space of differentials 
of $\cG$ over $\cH$, if $d$ is compatible with partial derivations and for any 
$\kappa[[e^\cH]]$-derivation $\delta\colon\kappa[[e^\cG]]\to V$ compatible with partial
derivations, there exists a unique $\kappa[[e^\cG]]$-linear map $f\colon\Omega_{\cG/\cH}\to V$
such that $f\circ d=\delta$. 
\end{defn}
In other words, the vector space $\Omega_{\cG/\cH}$ 
is the universal object in the category of $\kappa[[e^\cH]]$-derivations on
$\kappa[[e^\cG]]$ compatible with partial derivations. Elements of $\Omega_{\cG/\cH}$ are called 
differentials of $\cG$ over $\cH$, or simply differentials if $\cG$ and $\cH$ are obvious in the
context. They are different from K\"{a}hler differentials of
$\kappa[[e^\cG]]$ over $\kappa[[e^\cH]]$, which form the universal object
$\Omega_{\kappa[[e^\cG]]/\kappa[[e^\cH]]}$ in the category of
$\kappa[[e^\cH]]$-derivations on $\kappa[[e^\cG]]$. For example, we will see in the next
proposition that $\Omega_{{\mathbb Z}/0}$ is an one-dimensional $\kappa[[e^{\mathbb Z}]]$-vector
space with the usual order on $\mathbb Z$. However, if the characteristic of $\kappa$ is zero,
differential basis of $\Omega_{\kappa[[e^{\mathbb Z}]]/\kappa}$ (that is, a subset
$B$ of $\kappa[[e^{\mathbb Z}]]$ such that $\{d\Phi\colon\Phi\in B\}$ forms a basis of
$\Omega_{\kappa[[e^{\mathbb Z}]]/\kappa}$) is
exactly transcendence basis of $\kappa((X))$ over
$\kappa$, whose cardinality is infinite.  This example shows also that there do exist
$\kappa[[e^\cH]]$-derivations not compatible with partial derivations.

\begin{prop}[existence of differentials]
The vector space of differentials 
of $\cG$ over $\cH$ exists (with our assumption that $\cG/\cH$ is free of rank $n$). 
The differentials $d X_1,\cdots,d X_n$ form a basis of 
$\Omega_{\cG/\cH}$ for any set of variables $X_1,\cdots,X_n$.
\end{prop}
\begin{proof}
Let $V$ be a $\kappa[[e^\cG]]$-vector space with basis $v_1,\cdots,v_n$. The
$\kappa[[e^\cH]]$-derivation $d\colon\kappa[[e^\cG]]\to V$ defined by
$d\Phi=\sum(\partial\Phi/\partial X_i)v_i$ is compatible with 
$\partial\Phi/\partial X_1,\cdots,\partial\Phi/\partial X_n$ and hence compatible with
partial derivations. It is easy to see that $V$ together with $d$ 
satisfies the universal property.
\end{proof}
Let $\Phi$ be a positive generalized power series and $c_i$ be scalars. For the special
case that there are only finitely many non-zero $c_i$, clearly
$$
d(c_0+c_1\Phi+c_2\Phi^2+\cdots)=(c_1+2c_2\Phi+3c_3\Phi^2+\cdots)d\Phi.
$$
For arbitrary $c_i$, note that the $\kappa$-coefficients of the generalized power series
on both sides of the above equation at any $g\in\cG$ involve only finitely many $c_i$.
By reduction to the special case, we see that the above equation always holds. In particular 
$d\log(1+\Phi)=d\Phi/(1+\Phi)=d(1+\Phi)/(1+\Phi)$, 
if the characteristic of $\kappa$ is zero.
Even though logarithmic functions are not defined for a field with positive characteristic, we
still use the notation
$$
\dlog\Phi:=\frac{d\Phi}{\Phi}
$$
for a non-zero generalized power series $\Phi$ with coefficients in an arbitrary field. The
operator $\dlog$ transforms the multiplication of non-zero generalized power series to an addition:
$$
\dlog(\Phi_1\Phi_2)=\dlog\Phi_1+\dlog\Phi_2.
$$

Let $X_1,\cdots,X_n$ be variables of $\kappa[[e^\cG]]$ over $\kappa[[e^\cH]]$. Given
$\Phi_1,\cdots,\Phi_n\in\kappa[[e^\cG]]$, we define their Jacobian with respect to
$X_1,\cdots,X_n$ to be
$$
\left|\frac{\partial\bPhi}{\partial\bX}\right|:=
\left|\frac{\partial(\Phi_1,\cdots,\Phi_n)}{\partial(X_1,\cdots,X_n)}\right|:=
\det\left(\frac{\partial\Phi_i}{\partial X_j}\right).
$$
One is often interested in the $\kappa[[e^\cH]]$-coefficient of
\begin{equation}\label{eq:jac}
\frac{1}{\bPhi}\left|\frac{\partial\bPhi}{\partial\bX}\right|=
\frac{1}{\Phi_1\cdots \Phi_n}
\left|\frac{\partial(\Phi_1,\cdots,\Phi_n)}{\partial(X_1,\cdots,X_n)}\right|
\end{equation}
at $\bX^{-1}$ with the conventions $\bPhi:=\Phi_1\cdots \Phi_n$ and $\bX:=X_1\cdots X_n$. Since
a Jacobian appears in the generalized power series, it is more  natural to work on the $n$th
exterior product $\wedge^n\Omega_{\cG/\cH}$ of
$\Omega_{\cG/\cH}$, which is a dimension one $\kappa[[e^\cG]]$-vector space with a basis 
$$
d\bX:=d X_1\wedge\cdots\wedge d X_n.
$$ 
For $\Phi\in\kappa[[e^\cG]]$ and $g\in\cG$, the $\kappa$-coefficient of 
$\Phi d\bX$ at $e^gd\bX$ is defined as the $\kappa$-coefficient of $\Phi$ at
$e^g$; the $\kappa[[e^\cH]]$-coefficient
of $\Phi d\bX$ at $X_1^{j_1}\cdots X_n^{j_n}d\bX$ is defined as the 
$\kappa[[e^\cH]]$-coefficient of $\Phi$ at $X_1^{j_1}\cdots X_n^{j_n}$. 
\begin{prop}[vanishing of coefficients]\label{prop:coeff1}
Let $\Phi_i,\cdots,\Phi_n$ be non-zero generalized power series. If some $i_j$ 
is not equal to $1$, the $\kappa[[e^\cH]]$-coefficient of
$$
\frac{d\bPhi}{\bPhi^\bi}:=
\frac{d\Phi_1}{\Phi_1^{i_1}}\wedge\cdots\wedge\frac{d\Phi_n}{\Phi_n^{i_n}}
$$
at $\dlog\bX:=d\bX/\bX$ is zero for any set of variables $X_1,\cdots,X_n$.
\end{prop}
\begin{proof}
The proposition is equivalent to that the $\kappa$-coefficient $c_h$ of
$d\bPhi/\bPhi^\bi$ at $e^h\dlog\bX$ is zero for any $h\in\cH$. Since $c_h$ involves only
finitely many non-zero $\kappa$-coefficients of $\Phi_i$, we may assume that
$\Phi_1,\cdots,\Phi_n$ has only finitely many non-zero $\kappa$-coefficients $a_1,\cdots,a_m$.
The coefficient $c_h$ is obtained from $a_1,\cdots,a_m$ by finitely many algebraic
operations in $\kappa$ (additions, subtractions, multiplications and divisions). There is a
polynomial
$f\in{\mathbb Z}[T_1,\cdots,T_m]$ and $s\in\mathbb N$ such that 
$$
\frac{f(a_1,\cdots,a_m)}{(a_1\cdots a_m)^s}=c_h.
$$
To show $f(a_1,\cdots,a_m)$ is zero, it suffices to show that so is $f(T_1,\cdots,T_m)$.
Replacing $a_i$ by $T_i$, the $\kappa$-coefficients of $\Phi_i$ becomes elements in
the field ${\mathbb Q}(T_1,\cdots,T_m)$. So we may assume that $\kappa={\mathbb
Q}(T_1,\cdots,T_m)$. In particular, $\kappa$ has characteristic zero.

Following the idea of \cite[Section 1]{chen-mck-tow-wan-wri:rpserc}, we treat first the special
case that all $i_\ell$ are zero. The derivation $d$ is $\kappa$-linear, so we may assume
furthermore that $\Phi_i$ has only one non-zero $\kappa$-coefficient, that is,
$\Phi_i=a_ie^{h_i}X_1^{s_{i1}}\cdots X_n^{s_{in}}$ for some
$a_i\in\kappa$, $h_i\in\cH$ and $s_{ij}\in\mathbb Z$. Under these assumptions,
$$
d\Phi_1\wedge\cdots\wedge d\Phi_n=a_1\cdots a_n
(\det s_{ij})e^{\sum h_i}X_1^{-1+\sum s_{i1}}\cdots X_n^{-1+\sum s_{in}} d\bX.
$$
In order to have non-zero coefficients, $\det s_{ij}$ can not vanish in $\kappa$. But this
would imply that the power of some $X_i$ in the right hand side of the above equation is
not $-1$. Hence the $\kappa[[e^\cH]]$-coefficient of $d\bPhi/\bPhi^\bi$ at $\dlog\bX$ is
zero.

For the general case, we may assume that $i_\ell=1$ for $\ell$ greater than some fixed $j$
and $i_\ell\neq 1$ for $\ell\leq j$. Since the characteristic of $\kappa$
is assumed to be zero, $1-i_\ell$ is invertible in $\kappa$ for $\ell\leq j$. The general case
is reduced to the special case from the following straightforward computation.
\begin{eqnarray*}
&&
\frac{d\Phi_1}{\Phi_1^{i_1}}\wedge\cdots\wedge\frac{d\Phi_n}{\Phi_n^{i_n}}\\
&=&
d\left(\frac{\Phi_1^{1-i_1}}{1-i_1}\right)\wedge
\cdots\wedge d\left(\frac{\Phi_j^{1-i_j}}{1-i_j}\right)\wedge
\frac{d\Phi_{j+1}}{\Phi_{j+1}}\wedge\cdots\wedge\frac{d\Phi_n}{\Phi_n}\\
&=&
d\left(\frac{\Phi_1^{1-i_1}}{1-i_1}\Phi_{j+1}^{-1}\cdots\Phi_n^{-1}\right)\wedge
d\left(\frac{\Phi_2^{1-i_2}}{1-i_2}\right)\wedge
\cdots\wedge d\left(\frac{\Phi_j^{1-i_j}}{1-i_j}\right)\wedge\\
&&
d\Phi_{j+1}\wedge\cdots\wedge d\Phi_n.
\end{eqnarray*}
\end{proof}
\begin{prop}[determinant of exponents]\label{prop:coeff2}
Let $X_1,\cdots,X_n$ be a set of variables and $\Phi_1,\cdots,\Phi_n$ be non-zero generalized
power series with the factorizations
$\Phi_i=a_ie^{h_i}X_1^{s_{i1}}\cdots X_n^{s_{in}}(1+\tilde{\Phi}_i)$. 
The $\kappa[[e^\cH]]$-coefficient of 
$$
\dlog\bPhi:=
\dlog\Phi_1\wedge\cdots\wedge \dlog\Phi_n=
\frac{1}{\bPhi}\left|\frac{\partial\bPhi}{\partial\bX}\right|d\bX
$$
at $\dlog\bX$ is $\det s_{ij}$. 
\end{prop}
\begin{proof}
As the proof of Proposition~\ref{prop:coeff1}, we may assume that $\kappa$ has
characteristic zero. The element $\log(1+\tilde{\Phi}_i)$ can be defined, with which
$$
\dlog\Phi_i
=
\dlog(1+\tilde{\Phi}_i)+\sum_{j=1}^n s_{ij}\frac{d X_j}{X_j}.
$$
By Proposition~\ref{prop:coeff1}, the $\kappa[[e^\cH]]$-coefficient of $\dlog\bPhi$
at $\dlog\bX$ equals to that of
$$
\left(\sum_{j=1}^n s_{1j}\frac{d X_j}{X_j}\right)\wedge\cdots\wedge
\left(\sum_{j=1}^n s_{nj}\frac{d X_j}{X_j}\right),
$$
which is clearly $\det s_{ij}$.
\end{proof}
We would like to define a map independent of the choice of variables with the effect of
taking $\kappa[[e^\cH]]$-coefficients of a generalized power series. Restricting ourselves to
$\kappa[[e^\cG]]$ does not work. For instance, if we replace
$X_1$ by $X_1/\varphi$ with some non-zero element $\varphi\in\kappa[[e^\cH]]$, the
$\kappa[[e^\cH]]$-coefficient $\eta$ of a generalized power series at
$X_1^{j_1}\cdots X_n^{j_n}$ becomes $\eta\varphi^{j_1}$. Working on $\wedge^n\Omega_{\cG/\cH}$
still has problems: For instance, if we switch the order of two variables, the $\kappa[[e^\cH]]$-coefficients
of a generalized power series change signs. In particular, the $\kappa[[e^\cH]]$-coefficient
of $d X_1\wedge\cdots\wedge d X_n$ at 
$d X_2\wedge d X_1\wedge d X_3\wedge\cdots\wedge d X_n$ is $-1$. In the next section, we will
introduce parameters and generalized fractions to achieve our goal.


\section{Residues}\label{sec:pr}


Recall that $\cG$ is a totally ordered Abelian group generated freely by a subgroup 
$\cH$ and $n$ elements. 

\begin{defn}[multiplicity]
Let $\Phi$ be a non-zero generalized power series with the
factorization $\Phi=aY(1+\tilde{\Phi})$. The multiplicities of $\Phi$ with respect to a set of
variables $X_1,\cdots,X_n$ are the integers $i_1,\cdots,i_n$ such that 
$Y=e^hX_1^{i_1}\cdots X_n^{i_n}$ for some $h\in\cH$.
\end{defn}
\begin{defn}[parameter]\label{defn:logpar}
Non-zero generalized power series $\Phi_1,\cdots,\Phi_n$ form a
system of parameters (or simply parameters) of $\kappa[[e^\cG]]$ over $\kappa[[e^\cH]]$ if the
determinant of their multiplicities (with respect to a set of variables) is not zero in $\kappa$.
\end{defn}
The definition is independent of the choice of variables. The determinant of the multiplicities
of a system of parameters is not zero in 
$\mathbb Z$. Let $\Phi_1,\cdots,\Phi_n$ be parameters of $\kappa[[e^\cG]]$ over 
$\kappa[[e^\cH]]$ with the factorizations $\Phi_i=a_iY_i(1+\tilde{\Phi}_i)$. The 
necessary and sufficient condition for 
$e^{h_1}Y_1^{i_1}\cdots Y_n^{i_n}=e^{h_2}Y_1^{j_1}\cdots Y_n^{j_n}$ with $h_i\in\cH$ is
$i_1=j_1$, $\cdots$, $i_n=j_n$ and $h_1=h_2$.
\begin{lem}
Let $\Phi_1,\cdots,\Phi_n$ be parameters of $\kappa[[e^\cG]]$ over $\kappa[[e^\cH]]$.
Then $\dlog\bPhi$ is a basis of the $\kappa[[e^\cG]]$-vector space $\wedge^n\Omega_{\cG/\cH}$.
\end{lem}
\begin{proof}
Since $\wedge^n\Omega_{\cG/\cH}$ has dimension one, we need to check that $\dlog\bPhi\neq 0$.
Let $X_1,\cdots,X_n$ be variables. Replacing $X_i$ by its inverse if $\log X_i<0$, we may assume
that $\log X_i>0$ for all $i$. With the factorizations
$\Phi_i=a_ie^{h_i}X_1^{s_{i1}}\cdots X_n^{s_{in}}(1+\tilde{\Phi}_i)$, we can write
\begin{eqnarray*}
d\log\bPhi
&=&
(\frac{d\tilde{\Phi}_1}{1+\tilde{\Phi}_1}+\sum_{j}s_{1j}\frac{d X_j}{X_j})
\wedge\cdots\wedge
(\frac{d\tilde{\Phi}_n}{1+\tilde{\Phi}_n}+\sum_{j}s_{nj}\frac{d X_j}{X_j})\\
&=&                                                                                      
\Psi d\bX+\frac{\det s_{ij}}{\bX}d\bX
\end{eqnarray*}
for some $\Psi\in\kappa[[e^\cG]]$. Note that, from our convention of positivity
of $\tilde{\Phi}_i$, the support of $\Psi$ consists of only elements greater than
$-\log\bX$. Since
$\Phi_1,\cdots,\Phi_n$ are parameters, the leading coefficient
$\det s_{ij}$ of $\Psi+(\det s_{ij})\bX^{-1}$ is not zero. Therefore
$\dlog\bPhi\neq 0$.
\end{proof}

Let $V$ be a $\kappa[[e^\cG]]$-vector space. In the set
$$
\{(\alpha,\Phi_1,\cdots,\Phi_n)\in V\times\kappa[[e^\cG]]^n\colon
\text{ $\Phi_1,\cdots,\Phi_n$ are parameters}\},
$$
we define an equivalence relation:
$$
(\alpha,\Phi_1,\cdots,\Phi_n)\sim(\beta,\Psi_1,\cdots,\Psi_n)
\Longleftrightarrow
\frac{\beta}{\det t_{ij}}=\det u_{ij}\frac{\alpha}{\det s_{ij}},
$$
where $s_{ij}$ (resp. $t_{ij}$) are multiplicities of $\Phi_i$ (resp. $\Psi_i$) with respect to
a set of variables $X_1,\cdots,X_n$ (resp. $Y_1,\cdots,Y_n$)
and $Y_i=e^{h_i}X_1^{u_{i1}}\cdots X_n^{u_{in}}$. The equivalence relation is independent of the
choices of variables. 
\begin{defn}[generalized fraction]
A generalized fraction 
$$
\gfrac{\alpha}{\log\bPhi}:=\gfrac{\alpha}{\log\Phi_1,\cdots,\log\Phi_n}
$$
is the equivalence class containing $(\alpha,\Phi_1,\cdots,\Phi_n)$.
We call $\alpha$ the numerator of the generalized fraction. The set of generalized fractions with
numerators in $V$ is denoted by $\mH(V)$. 
\end{defn}
We choose the notation $\mH(V)$, because it might relate to
some cohomology object as the case of the theory of local cohomology residues for formal power
series rings.

\begin{defn}[residue]
Let $X_1,\cdots,X_n$ be variables of $\kappa[[e^\cG]]$ over $\kappa[[e^\cH]]$.
We define the residue map
$$
\res_{X_1,\cdots,X_n}\colon\mH(\wedge^n\Omega_{\cG/\cH})\to\kappa[[e^\cH]]
$$ 
with respect to $X_1,\cdots,X_n$ by
$$
\res_{X_1,\cdots,X_n}\gfrac{\Phi\dlog\bX}{\log\bX}=
\text{ the constant term of $\Phi$ in $\kappa[[e^\cH]]$,}
$$
where $\Phi\in\kappa[[e^\cG]]$.
\end{defn}
\begin{prop}[invariance of residues]
$\res_{X_1,\cdots,X_n}=\res_{Y_1,\cdots,Y_n}$ for any two sets of variables $X_1,\cdots,X_n$
and $Y_1,\cdots,Y_n$ of $\kappa[[e^\cG]]$ over $\kappa[[e^\cH]]$.
\end{prop}
\begin{proof}
Write $Y_i=e^{h_i}X_1^{s_{i1}}\cdots X_n^{s_{in}}$. Then
$$
\dlog\bY=
\dlog(X_1^{s_{11}}\cdots X_n^{s_{1n}})\wedge\cdots\wedge\dlog(X_1^{s_{n1}}\cdots X_n^{s_{nn}})
=(\det s_{ij})\dlog\bX.
$$
For any $\Phi\in\kappa[[e^\cG]]$,
$$
\gfrac{\Phi\dlog\bY}{\log\bY}=\gfrac{(\det s_{ij})^{-1}\Phi\dlog\bY}{\log\bX}
=\gfrac{\Phi\dlog\bX}{\log\bX}.
$$
Therefore $\res_{X_1,\cdots,X_n}=\res_{Y_1,\cdots,Y_n}$.
\end{proof}
Taking the residue is a map equating coefficients independent of the choice of variables. We denote
$\res:=\res_{X_1,\cdots,X_n}$.

Let $\Phi_1,\cdots,\Phi_n\in\kappa[[e^\cG]]$ and $\varphi_{i_1,\cdots,i_n}\in\kappa[[e^\cH]]$, 
where the indices $i_\ell\in\mathbb Z$. We assume that there are only
finitely many $\varphi_{i_1,\cdots,i_n}\Phi_1^{i_1}\cdots\Phi_n^{i_n}$ whose support
contains any fixed $g\in\cG$. Under such an assumption, we can define an element
$$
\Psi(\bPhi):=\sum\varphi_{i_1,\cdots,i_n}\Phi_1^{i_1}\cdots\Phi_n^{i_n}\in\kappa\{e^\cG\},
$$
whose $\kappa$-coefficient at $e^g$ is the sum of the $\kappa$-coefficients of all
$\varphi_{i_1,\cdots,i_n}\Phi_1^{i_1}\cdots\Phi_n^{i_n}$ at $e^g$. We say that
$\Psi(\bPhi)$ is represented by $\Phi_1,\cdots,\Phi_n$ with $\kappa[[e^\cH]]$-coefficient
$\varphi_{i_1,\cdots,i_n}$ at $\Phi_1^{i_1}\cdots\Phi_n^{i_n}$. Given
$\varphi_{i_1,\cdots,i_n}$, whether or not $\Psi(\bPhi)\in\kappa[[e^\cG]]$ depends on the
sequence $\Phi_1,\cdots,\Phi_n$. For instance, as seen in  Example~\ref{example:z2}, the 
inverse of $X+Y$ in the iterated Laurent series $\kappa[[e^{\cG_1}]]$ is
$$
\Psi(X,Y):=Y^{-1}-XY^{-2}+X^2Y^{-3}-\cdots.
$$
However,
$$
\Psi(Y,X)=X^{-1}-YX^{-2}+Y^2X^{-3}-\cdots
$$
although defined is not contained in $\kappa[[e^{\cG_1}]]$.

If there are only finitely many nonzero $\kappa[[e^\cH]]$-coefficients for a representation
$\Psi(\bPhi)$ of an element in $\kappa\{e^\cG\}$ and the indices of these
nonzero coefficients are all non-negative, we say that
$\Psi(\bPhi)$ is a polynomial in $\Phi_1,\cdots,\Phi_n$ with coefficients in $\kappa[[e^\cH]]$.
The set of these polynomials is a subring of $\kappa[[e^\cG]]$ denoted by
$\kappa[[e^\cH]][\Phi_1,\cdots,\Phi_n]$, which is exactly the image of the homomorphism
$$
\kappa[[e^{\cH}]][Y_1,\cdots,Y_n]\to\kappa[[e^\cG]]
$$
of $\kappa[[e^\cH]]$-algebras sending $Y_i$ to $\Phi_i$. For $F\in\kappa[[e^{\cH}]][Y_1,\cdots,Y_n]$
sending to $\Psi(\bPhi)$ under this homomorphism, we use also the notation $F(\bPhi):=\Psi(\bPhi)$.
If we assume furthermore that
$\Phi_1,\cdots,\Phi_n$ are parameters, then the above homomorphism is ono-to-one. A representation
$\Psi(\bPhi)$ of an element in $\kappa[[e^\cG]]$ is a rational function in
$\Phi_1,\cdots,\Phi_n$ with coefficients in $\kappa[[e^\cH]]$ if there exist 
$F_1,F_2\in\kappa[[e^{\cH}]][Y_1,\cdots,Y_n]$ with $F_1\neq 0$ such that $F_1(\bPhi)\Psi(\bPhi)=F_2(\bPhi)$.
Let $F=F_1/F_2\in\kappa[[e^{\cH}]](Y_1,\cdots,Y_n)$.
We use also the notation $F(\bPhi):=\Psi(\bPhi)$. 

Now we interpret and generalize Jacobi's formula.

\begin{thm}[Jacobi]\label{thm:Jacobi} 
Given a representation
$\Psi(\bPhi)=\sum\varphi_{i_1,\cdots,i_n}\Phi_1^{i_1}\cdots\Phi_n^{i_n}\in\kappa[[e^\cG]]$
by parameters $\Phi_1,\cdots,\Phi_n$ with $\varphi_{i_1,\cdots,i_n}\in\kappa[[e^\cH]]$,
$$
\res\gfrac{\Psi(\bPhi)\dlog\bPhi}{\log\bPhi}=\varphi_{0,\cdots,0}.
$$
\end{thm}
\begin{proof}
Let $s_{ij}$ be the multiplicities of $\bPhi$ with respect to a set of variables $\bX$. Since
$$
\gfrac{\Psi(\bPhi)\dlog\bPhi}{\log\bPhi}
=
\gfrac{\frac{\Psi(\bPhi)}{\bPhi}
\left|\frac{\partial\bPhi}{\partial\bX}\right|d\bX}{\log\bPhi}
=
\gfrac{\frac{\bX}{\det s_{ij}}
\frac{\Psi(\bPhi)}{\bPhi}
\left|\frac{\partial\bPhi}{\partial\bX}\right|\dlog\bX}{\log\bX},
$$
the theorem is equivalent to the claim that the $\kappa$-coefficients of
$(\det s_{ij})\varphi_{0,\cdots,0}$ at $e^h$ and
$(\Psi(\bPhi)/\bPhi)|\partial\bPhi/\partial\bX|$ 
at $e^h\bX^{-1}$ are the same for any $h\in\cH$. The latter involves
only finitely many $\varphi_{i_1,\cdots,i_n}$, so we may assume that only
finitely many $\varphi_{i_1,\cdots,i_n}$ are not zero. From linearity, we may assume furthermore
that
$\Psi=X_1^{i_1}\cdots X_n^{i_n}$. The theorem in such a special case was proved in
Propositions~\ref{prop:coeff1} and \ref{prop:coeff2}.
\end{proof}

Another interpretation and generalization of Jacobi's formula in characteristic zero
can be found in \cite[Theorem 3.7]{xin:rtmns}. While \cite{xin:rtmns} investigates the interplay
of two fields, we work on one vector space of differentials. In our
approach, combinatorial information appears naturally through a residue map with Jacobians 
resulted from parameters changes.
In \cite{xin:rtmns}, Jacobi's formula is called a residue theorem. However,
residue theorem usually refers to Cauchy's theorem, which counts residues of a 
meromorphic function in a region. As a global result relating the poles of a 
meromorphic function, Cauchy's residue theorem is considered in a very general
context by Grothendieck in algebraic geometry. Jacobi's formula, exploring parameters changes
of  one point, is merely a local property!

Formulae of the Lagrange inversion type can be studied in the field of generalized power series.
Along this direction, one needs to know whether or not every generalized power series
can be represented by a system of parameters and in what sense a representation is unique. 
Let $\Phi_1,\cdots,\Phi_n$ be a system of parameters of $\kappa[[e^\cG]]$ over $\kappa[[e^\cH]]$ with the
factorizations $\Phi_i=a_iY_i(1+\tilde{\Phi}_i)$. The following uniqueness property is obvious: If 
$$
\sum\varphi^{(1)}_{i_1,\cdots,i_n}Y_1^{i_1}\cdots Y_n^{i_n}
=\sum\varphi^{(2)}_{i_1,\cdots,i_n}Y_1^{i_1}\cdots Y_n^{i_n}\in\kappa[[e^\cG]],
$$
then $\varphi^{(1)}_{i_1,\cdots,i_n}=\varphi^{(2)}_{i_1,\cdots,i_n}$ for all $i_1,\cdots,i_n$.

\begin{defn}[regular parameter]
A system of parameters $\Phi_1,\cdots,\Phi_n$ of $\kappa[[e^\cG]]$ over $\kappa[[e^\cH]]$ with the
factorizations $\Phi_i=a_iY_i(1+\tilde{\Phi}_i)$ is regular if, for every element $\Psi\in\kappa[[e^\cG]]$,
there exists an unique element $\sum\varphi_{i_1,\cdots,i_n}Y_1^{i_1}\cdots Y_n^{i_n}\in\kappa[[e^\cG]]$
such that $\Psi=\sum\varphi_{i_1,\cdots,i_n}\Phi_1^{i_1}\cdots\Phi_n^{i_n}$.
\end{defn}

Clearly, variables are regular parameters. 

\begin{prop}[characterization of regularity]
A system of parameters of $\kappa[[e^\cG]]$ over $\kappa[[e^\cH]]$ is regular if and only if the
determinant of their multiplicities (with respect to a set of variables) is invertible in
$\mathbb Z$.
\end{prop}  
\begin{proof}
Let $\Phi_1,\cdots,\Phi_n$ be parameters of $\kappa[[e^\cG]]$ over $\kappa[[e^\cH]]$ with the
factorizations $\Phi_i=a_iY_i(1+\tilde{\Phi}_i)$. 

Assume that the determinant of their
multiplicities is invertible in $\mathbb Z$. This assumption is equivalent to that
$Y_1,\cdots,Y_n$ are variables. We need to find the $\kappa$-coefficient $a_g$ of 
a generalized power series $\sum\varphi_{i_1,\cdots,i_n}Y_1^{i_1}\cdots Y_n^{i_n}$ at $e^g$
for each $g\in\cG$ to represent a given generalized power series $\Psi=\sum b_ge^g$. Let 
$$
A=\supp\tilde{\Phi}_1+\cdots+\supp\tilde{\Phi}_n,
$$
$iA$ be the sum of $i$ copies of $A$ for $i>0$, and $0A:=\{0\}$. The well-ordered
set $\bar{A}:=\cup_{i\geq 0}iA$ contains $\supp(1+\tilde{\Phi}_1)^{i_1}\cdots(1+\tilde{\Phi}_n)^{i_n}$.
Let $B$ be a well-ordered set containing $\supp\Psi$. If $g\not\in\bar{A}+B$, we define 
$a_g=0$. For $g\in\bar{A}+B$, we consider the equation $x+y=g$ subject to the 
constraints $x\in\bar{A}+B$ and $y\in\supp(1+\tilde{\Phi}_1)^{i_1}\cdots(1+\tilde{\Phi}_n)^{i_n}$, 
where $i_1,\cdots,i_n$ are integers satisfy $x=h+i_1\log Y_1+\cdots+i_n\log Y_n$ for
some $h\in\cH$. By Lemma~\ref{lem:9413451}, the equation with the constraints
has finitely many solutions. If $(x,y)=(g,0)$ is the only solution, for instance if $g$ is the
smallest element of $\bar{A}+B$, we define 
$a_g:=b_g$. If it has other solutions, say 
$(g_{11},g_{21})$, ... , $(g_{1m},g_{2m})$ besides $(g,0)$, we would like to define
$$
a_g:=b_g-\sum_{\ell=1}^ma_{g_{1\ell}}c_\ell
$$
inductively in terms of $a_{g_{1\ell}}$, where $g_{1\ell}=h_\ell+i_{\ell 1}\log Y_1+\cdots+i_{\ell n}\log Y_n$,
$h_\ell\in\cH$, $i_{\ell 1},\cdots,i_{\ell n}\in\mathbb Z$ and $c_\ell$ is the $\kappa$-coefficient of
$(1+\tilde{\Phi}_1)^{i_{\ell 1}}\cdots(1+\tilde{\Phi}_n)^{i_{\ell n}}$ at $e^{g_{2\ell}}$. To see
the inductive process working, we observe that $g_{1\ell}<g$, since $g_{2\ell}>0$. 
Moreover, if the equation $x+y=g_{1\ell}^{(1)}$ for each $g_{1\ell}^{(1)}:=g_{1\ell}$
subject to the constraints above has only one solution, $a_{g_{1\ell}}$
has been defined. Let $(g_{11}^{(2)},g_{21}^{(2)})$, $(g_{12}^{(2)},g_{22}^{(2)})$,
$(g_{13}^{(2)},g_{23}^{(2)})$ ...  be solutions of other equations if any. 
We repeat the process for the equations $x+y=g_{1\ell}^{(2)}$ with the same constraints. 
If these equations have more than one solution, we continue the process. The process has 
to stopped in finitely many steps, since the elements $g_{1\ell}^{(i)}$ obtained are contained 
in $\bar{A}+B$, which consists of no strictly decreasing infinite sequences. Therefore $a_g$ is defined.
From the construction, $\Psi$ is represented by $\Phi_1,\cdots,\Phi_n$ with 
$\varphi_{i_1,\cdots,i_n}$, where $\sum a_ge^g=\sum\varphi_{i_1,\cdots,i_n}Y_1^{i_1}\cdots
Y_n^{i_n}\in\kappa[[e^\cG]]$. 

Assume that there are two representations
$$
\Psi=\sum\varphi^{(1)}_{i_1,\cdots,i_n}\Phi_1^{i_1}\cdots\Phi_n^{i_n}=
\sum\varphi^{(2)}_{i_1,\cdots,i_n}\Phi_1^{i_1}\cdots\Phi_n^{i_n},
$$
with $\Psi_i=\sum\varphi^{(i)}_{i_1,\cdots,i_n}Y_1^{i_1}\cdots Y_n^{i_n}\in\kappa[[e^\cG]]$.
In the above process, we may take $B=\supp\Psi\cup\supp\Psi_1\cup\supp\Psi_2$, which contains
both $\supp\Psi_1$ and $\supp\Psi_2$. As
$a_g$ is determined by $b_g$ for $g\in\bar{A}+B$, the representations must be the same.

Now we assume that $\Phi_1,\cdots,\Phi_n$ are regular parameters. Let $X_1,\cdots,X_n$ be a set
of variables. There exist $\Psi_\ell=\sum\varphi^{(\ell)}_{i_1,\cdots,i_n}Y_1^{i_1}\cdots
Y_n^{i_n}\in\kappa[[e^\cG]]$ such that 
$X_\ell=\sum\varphi^{(\ell)}_{i_1,\cdots,i_n}\Phi_1^{i_1}\cdots\Phi_n^{i_n}$.
Since $Y_1,\cdots,Y_n$ are parameters, the minimal element of 
$\supp\Psi_\ell$ is $\log X_\ell$. This implies that 
$X_\ell=e^{h_\ell}Y_1^{i_{\ell 1}}\cdots Y_n^{i_{\ell n}}$
for some $h_\ell\in\cH$ and $i_{\ell 1},\cdots,i_{\ell n}\in\mathbb Z$. Therefore the
determinant of the multiplicities of $\Phi_1,\cdots,\Phi_n$ is invertible in $\mathbb Z$.
\end{proof}

The theme of Lagrange inversions in the context of generalized power series is the
interrelations between two systems of regular parameters. Let $\Phi_1,\cdots,\Phi_n$ be a
system of regular parameters represented by another system of regular parameters
$\Upsilon_1,\cdots,\Upsilon_n$. The expression
$$
\res\gfrac{\frac{\Upsilon_\ell}{\Phi_1^{i_1}\cdots\Phi_n^{i_n}}\dlog\bPhi}{\log\bPhi}
$$
gives the $\kappa[[e^\cH]]$-coefficient of $\Upsilon_\ell$ at $\Phi_1^{i_1}\cdots\Phi_n^{i_n}$.
Properties of generalized fraction and residues can be used to compute the coefficient in terms
of $\kappa[[e^\cH]]$-coefficients of $\Phi_1,\cdots,\Phi_n$ at monomials in
$\Upsilon_1,\cdots,\Upsilon_n$.


\section{Dyson's conjecture}\label{sec:appl}


Let $a_1,\cdots,a_n$ be non-negative integers. Dyson's conjecture \cite{dys:stelcs} that 
$$
\text{the constant term of }\prod_{1\leq i\neq j\leq n}(1-\frac{X_i}{X_j})^{a_i}
=\frac{(a_1+\cdots+a_n)!}{a_1!\cdots a_n!}
$$
was confirmed by Wilson~\cite{wil:pcd} and Gunson~\cite{gun:pcdstel} independently of each other.
We interpret two known proofs of the Dyson's conjecture in terms of generalized power series
with coefficients in $\mathbb Q$ and exponents in ${\mathbb Z}^n$, which has a total order
compatible with the group structure (for instance, the lexicographic order). Let
$X_1,\cdots,X_n$ be variables of ${\mathbb Q}[[e^{{\mathbb Z}^n}]]$ over $\mathbb Q={\mathbb
Q}[[e^0]]$. 

In the first proof, we assume that the variables satisfy $\log X_1>\dots>\log X_n$.
Let $\Phi_i=\prod_{j=1,j\neq i}^{n}(1-X_i/X_j)^{-1}$. Using Lagrange interpolation, one
can show $\sum_{i=1}^n\Phi_i=1$. Wilson's proof to the Dyson's conjecture is
based on the parameters $X_1,\Phi_2,\cdots,\Phi_n$, whose multiplicities with respect to 
$X_1,\cdots,X_n$ have determinant 
$$
\det\left(\begin{matrix}
1       & 0     & 0     & \cdots & 0\\
1       & -1    & 0     & \cdots & 0\\
1       & 1     & -2     & \cdots & 0\\
\vdots  & \vdots & \vdots & \ddots & \vdots\\
1       & 1     & 1     & \cdots & -(n-1)
\end{matrix}\right)
=(n-1)!(-1)^{n-1}.
$$ 
Wilson's computation~\cite[Proof of Lemma 3]{wil:pcd} carried over to our context 
shows
$$
\dlog X_1\wedge\dlog\Phi_2\wedge\cdots\wedge\dlog\Phi_n
=c(n-1)!(-1)^{n-1}\Phi_1\dlog\bX
$$
for some $c\in\kappa$. The scalar $c$ is not zero, since 
$\dlog X_1\wedge\dlog\Phi_2\wedge\cdots\wedge\dlog\Phi_n$ 
generates $\wedge^n\Omega_{\cG/\cH}$. 
Let $\Psi(X_1,\cdots,X_n)=X_2^{-a_2}\cdots
X_n^{-a_n}(1-X_2-\cdots-X_n)^{-a_1}$. What we need to compute is the constant
term of $\Psi(X_1,\Phi_2,\cdots,\Phi_n)$, that is, the residue of
\begin{eqnarray*}
&&
\gfrac{\Psi(X_1,\Phi_2,\cdots,\Phi_n)\dlog\bX}{\log\bX}\\ 
&=&
\gfrac{\frac{c^{-1}\Psi(X_1,\Phi_2,\cdots,\Phi_n)}{1-\Phi_2-\cdots-\Phi_n}
\dlog X_1\wedge\dlog\Phi_2\wedge\cdots\wedge\dlog\Phi_n}
{\log X_1,\log\Phi_2,\cdots,\log\Phi_n}.
\end{eqnarray*}
By Theorem~\ref{thm:Jacobi}, the constant term of $\Psi(X_1,\Phi_2,\cdots,\Phi_n)$
is the same as that of
$$
\frac{c^{-1}\Psi(X_1,\cdots,X_n)}{1-X_2-\cdots-X_n}=
\frac{c^{-1}}{X_2^{a_2}\cdots X_n^{a_n}}
\sum_{k=0}^\infty\binom{k+a_1}{a_1}(X_2+\cdots+X_n)^k,
$$
which is
$$
c^{-1}\binom{a_1+\cdots+a_n}{a_1}
\binom{a_2+\cdots+a_n}{a_2,\cdots,a_n}\\
=
c^{-1}\frac{(a_1+\cdots+a_n)!}{a_1!\cdots a_n!}
$$
occurring when $k=a_2+\cdots+a_n$. Now Dyson's conjecture for the trivial case
$a_0=\cdots=a_n=0$ shows $c=1$.

In the second proof, we assume that $\log X_1<\dots<\log X_n$. Following Egorychev
\cite{eg:irccs}, we use the parameters 
$\Upsilon_i=(-1)^{i-1}X_i^{n-1}\prod_{j<k,\,\,j\neq i,\,\,k\neq i}(X_j-X_k)$, 
whose multiplicities with respect to $X_1,\cdots,X_n$ have determinant
$$
\det\left(\begin{matrix}
n-1    & n-2    & n-3     & \cdots & 0\\
n-2    & n-1    & n-3     & \cdots & 0\\
n-2    & n-3    & n-1     & \cdots & 0\\
\vdots & \vdots & \vdots & \ddots & \vdots\\
n-2    & n-3    & n-4     & \cdots & n-1
\end{matrix}\right)
=\frac{n!(n-1)}{2}
$$ 
\cite[proof of Theorem 5.3]{xin:rtmns}.
(The matrix has diagonal entries $n-1$ and other entries in each row, except the diagonal, are 
$n-2,n-3,\cdots,0$ from left to right.)  By Cramer's rule,
$$
\frac{\Upsilon_1}{X_1^i}+\cdots+\frac{\Upsilon_n}{X_n^i}=
\begin{cases}
\Delta:=\prod_{j<k}(X_j-X_k), & \text{ if $i=0$;}\\
0,      & \text{ if $i=1,\cdots,n-1$.}
\end{cases}
$$
Applying the derivation, we obtain
\begin{eqnarray*}
\Upsilon_1\dlog\Upsilon_1+\cdots+\Upsilon_n\dlog\Upsilon_n&=&
X_1\frac{\partial\Delta}{\partial X_1}\dlog
X_1+\cdots+X_n\frac{\partial\Delta}{\partial X_n}\dlog X_n,\\
\frac{\Upsilon_1}{X_1^i}\dlog\Upsilon_1+\cdots+\frac{\Upsilon_n}{X_n^i}\dlog\Upsilon_n&=&
i\frac{\Upsilon_1}{X_1^i}\dlog X_1+\cdots+i\frac{\Upsilon_n}{X_n^i}\dlog X_n
\end{eqnarray*}
($i=1,\cdots,n-1$). Exterior products of the above elements give rise to
$$
\frac{\Upsilon_1\cdots\Upsilon_n\Delta}{(X_1\cdots X_n)^{n-1}}\dlog\bUpsilon=
\frac{(n-1)!\Upsilon_1\cdots\Upsilon_n}{(X_1\cdots X_n)^{n-1}}
(X_1\frac{\partial\Delta}{\partial X_1}+\cdots+X_n\frac{\partial\Delta}{\partial X_n})\dlog\bX.
$$
Since
$$
X_1\frac{\partial\Delta}{\partial X_1}+\cdots+X_n\frac{\partial\Delta}{\partial X_n}=
\binom{n}{2}\Delta,
$$
the combinatorial number $n!(n-1)/2$ is exactly compensated in the identity
$$
\dlog\bUpsilon=\frac{n!(n-1)}{2}\dlog\bX.
$$
As \cite{xin:rtmns}, we condiser $\Psi(X_1,\cdots,X_n)=(X_1+\cdots+X_n)^{a_1+\cdots+a_n}/(X_1^{a_1}\cdots X_n^{a_n})$.
Since
$$
\prod_{i=1,i\neq j}^{n}(1-\frac{X_i}{X_j})=\frac{\Upsilon_1+\cdots+\Upsilon_n}{\Upsilon_j},
$$
we need to show
$$
\res\gfrac{\Psi(\bUpsilon)\dlog\bX}{\log\bX}=\frac{(a_1+\cdots+a_n)!}{a_1!\cdots a_n!}.
$$
This is a special case of Theorem~\ref{thm:Jacobi}, since
$$ 
\gfrac{\Psi(\bUpsilon)\dlog\bX}{\log\bX}=
\gfrac{\Psi(\bUpsilon)\dlog\bUpsilon}{\log\bUpsilon}
$$
and the constant term of $\Psi(X_1,\cdots,X_n)$ is $(a_1+\cdots+a_n)!/(a_1!\cdots a_n!)$.

One more proof by local cohomology residues is available. See
\cite[Identity 14]{hu:rmca}.



\end{document}